\documentclass[11pt]{article}
\usepackage{amssymb,amscd,latexsym}
\usepackage[mathscr]{euscript}

\DeclareFontFamily{OT1}{rsfs}{}
\DeclareFontShape{OT1}{rsfs}{n}{it}{<-> rsfs10}{}
\DeclareMathAlphabet{\curly}{OT1}{rsfs}{n}{it}

\newcommand\C{\mathbb C\,}
\newcommand\Pee{\mathbb P}
\newcommand\OO{\mathscr O}
\newcommand\I{\curly I}
\newcommand\AND{\!&\!}
\newcommand\into{\hookrightarrow}
\newcommand\so{\ \ \Longrightarrow\ }
\newcommand\res{\arrowvert_}
\newcommand\To{\longrightarrow}
\newcommand\End{\mathrm{End\,}}

\makeatletter \@addtoreset{equation}{section} \makeatother

\newtheorem{Theorem}[equation]{Theorem}
\newtheorem{Lemma}[equation]{Lemma}

\newenvironment{Proof}{\noindent\emph{Proof.}}{\hfill$\square$\\}
\newenvironment{Remarks}{\noindent\textbf{Remarks}.}{\\}

\title{An obstructed bundle on a Calabi-Yau 3-fold}
\author{R.\,P. Thomas} \date{}

\begin{document}
\maketitle

\begin{abstract} \noindent
Mirror symmetry suggests that on a Calabi-Yau 3-fold moduli
spaces of stable bundles, especially those with
degree zero and indivisible Chern class, might be smooth
(i.e. unobstructed, though perhaps of too high a dimension). This is
because smoothly embedded special Lagrangian cycles in the mirror
have unobstructed deformations. As there does not seem
to be a counterexample in the literature we provide one
here, showing that such a Tian-Todorov-McLean-type result cannot
hold.
\end{abstract}

\section{Introduction}

The existence of curves $C$ in Calabi-Yau 3-folds $X$ whose
deformations are obstructed has been known for some time.  This means
that the curve deforms to first order -- the linearisation of the
deformation problem has an $H^0(\nu)$ of solutions (where $\nu$ is
the normal bundle to $C\subset X$) -- but the actual deformation problem
has no solutions (in some direction in $H^0(\nu)$ at least). Thus the
moduli space of curves is not smooth at this point but has nilpotent
directions.

This implies that moduli spaces of sheaves can also be obstructed
(by taking the ideal sheaves
of the curves). However, the now folklore philosophy (see for instance
[Ty, Va]) of the relationship between the Kontsevich and
Strominger-Yau-Zaslow approaches to mirror symmetry suggests that for
bundles the situation might be different. Namely, under certain
circumstances, degree 0 \emph{stable} (Hermitian-Yang-Mills)
bundles on a Calabi-Yau 3-fold $X$ should correspond to
special Lagrangian cycles (with flat line bundles on them) on the
mirror $\widetilde X$. This would be set up by some sort of
analytical version of
the relative Fourier-Mukai transform for the dual $T^3$-fibrations of
$X$ and $\widetilde X$ whose existence is conjectured by
Strominger-Yau-Zaslow. And
McLean [Mc] has shown that deformations of smoothly embedded special
Lagrangians are unobstructed.

In general the picture is more complicated, with complexes of sheaves
on one side corresponding to elements of the derived category of the Fukaya
category of all Lagrangians on the other. Sheaves on one side can be
obstructed, while on the other singular Lagrangians and
``multiply wrapped'' cycles (images of degree\,$>1$ maps, or cycles with
higher rank vector bundles on them) may not have smooth moduli. But
the question remains as to what the analogue of the McLean's result is
in terms of sheaves on $X$. A natural subcategory of the derived
category of coherent sheaves on $X$ is the category of holomorphic
vector bundles, and while the mirror situation to having
multi-wrapped components is hard to understand, it is clearly related
to divisibility of the homology classes concerned. Thus it was
conjectured that moduli spaces of stable holomorphic bundles of degree
zero and indivisible Chern class on Calabi-Yau 3-folds might be
smooth. This would be analogous to McLean's result and also the
Tian-Todorov theorem that moduli spaces of the Calabi-Yau manifolds
themselves are smooth, again despite the fact that they can be of the
wrong dimension with a non-zero obstruction space.

Here, however, we will show that this cannot be true. We work on a
very particular smooth (3,3)-divisor $X\subset P=\Pee^2\times\Pee^2$
described below.
By the Lefschetz hyperplane theorem [GH], $H^2(X)\cong H^2(P)$ is
generated, over the integers, by the standard K\"ahler forms
$\omega_1,\ \omega_2$ pulled back from the two $\Pee^2$ factors of $P$.
By Poincar\'e duality $H^4(X)$ is generated over the integers by
$\omega_1^2/3$ and $\omega_2^2/3$, with $\omega_1\omega_2=\omega_1^2
+\omega_2^2$.

\begin{Theorem} \label{thm}
There exists a rank 2 holomorphic vector bundle $A$ on $X$ such that
\begin{itemize} \item $c_1(A)=0$, $c_2(A)=\omega_1^2+4/3\,\omega_2^2$
is indivisible,
\item $A$ is slope-stable (and so Gieseker-stable) with respect to the
K\"ahler forms $N\omega_1+\omega_2,\ N>1$, and
\item deformations $H^1(\End A)\cong\C^2$ are completely obstructed
to second order. \end{itemize}
Thus the component of the moduli space that
$A$ sits in is a thickened point. 
\end{Theorem}

In the published version of this paper, it was incorrectly claimed that this implied that the moduli space is isomorphic to $\mathrm{Spec}\,\C[\epsilon,\eta]/(\epsilon^2,\epsilon\eta,\eta^2)$. As Yukinobu Toda pointed out, this is not the critical locus of an algebraic function (as we now know all moduli of sheaves on Calabi-Yau 3-folds are). At the end of the paper we verify that in fact the moduli space is isomorphic to
$$
\mathrm{Spec}\,\C[\epsilon,\eta]/(\epsilon^2,\eta^2)\ =\ \mathrm{Crit}\,(\epsilon^3+\eta^3).
$$

\smallskip
The strategy to prove Theorem \ref{thm} is to employ the Serre construction relating vector
bundles to curves via zeros of a section. The work is in finding an
obstructed (rational) curve $C$ such that a bundle $A$ can be
constructed from it, and such that $A$ is stable (rational curves
usually correspond to unstable bundles)
with deformation theory exactly that of $C$.

It would be interesting to know what the mirror object to $A$
is. Although $c_2(A)$ is indivisible, it is the sum of two divisible
classes which are positive in some sense, and so may correspond to
some kind of sum of two multiply wrapped cycles, or some more
complicated (obstructed) object of the derived category of the
Fukaya category of the mirror. Also it is natural to expect the
correspondence between stable bundles and special Lagrangians to hold
only in the adiabatic limit of small SYZ $T^3$-fibres, where it is
possible that $A$ is not stable.\\

\noindent \textbf{Acknowledgements.}  I would like to thank Sheldon
Katz for pointing out [K] which prompted this work (though its
influence is well hidden now), and the organisers of the Harvard
Winter School on Mirror Symmetry in January 1999 for the stimulating
atmosphere in which the prompting occured. Thanks also to Ivan Smith
for useful conversations, and to Yukinobu Toda for pointing out an error
in the published version.

\section{An obstructed rational curve}

Consider a Calabi-Yau 3-fold $X=X_{3,3}\subset P$, a $(3,3)$ divisor
in $P=\Pee^2\times\Pee^2=\Pee^2_1\times\Pee^2_2$. This is elliptically
fibred over the second factor $\Pee^2_2$, and we will construct a
smooth example with a degenerate fibre $f$ that is the union of a line
with a double (thickened) line
$$
f=\{x^2y=0\}\subset\Pee^2_1.
$$
Moving away from this fibre we will have either smooth tori fibres, or
degenerate fibres that are unions of smooth conics and lines -- thus
the line $C=\{x=0\}$ whose double lies in $f$ cannot
deform inside $X$ beyond its first order deformations inside the
thickening.

Fixing local coordinates $(u,v)$ on the base $\Pee^2_2$, and
coordinates $[x:y:z]$ on $\Pee^2_1$, a local example of such a $(3,3)$
divisor is given by the zero locus of
$$
F=(x^2+z^2v)y+z^3u+(x+y)^3v+p,
$$
where the degenerate fibre $f$ is at $u=0=v$, and $p$ is any small
perturbation \emph{vanishing on $f$}. We have
\begin{eqnarray} \nonumber
\left.{\partial F\over\partial u}\right\res{u=0=v}\AND=\AND z^3+
{\partial p\over\partial u}, \\ \nonumber
\left.{\partial F\over\partial v}\right\res{u=0=v}\AND=\AND z^2y+
(x+y)^3+{\partial p\over\partial v}.
\end{eqnarray}
Thus for $p$ sufficiently small that $\partial p/\partial u$ is
pointwise smaller in norm than $z^3$ away from a fixed small
neighbourhood of $\{z=0\}\subset\Pee^2_1$, and $\partial p/\partial v$
pointwise smaller in norm than $z^2y+(x+y)^3$ away from a fixed small
neighbourhood of $\{z^2y+(x+y)^3=0\}$, we have $X=F^{-1}(0)$ smooth in
a neighbourhood of $u=0=v$\,: $\partial F/\partial u$ and $\partial
F/\partial v$ vanish simultaneously only in a small neighbourhood of the
intersection of $z=0$ and $z^2y+(x+y)^3=0$, but this is $[1:-1:0]\in
\Pee^2_1\times\{(0,0)\}$ which does not lie on $X$.

Globally then, with respect to coordinates $[x:y:z]$ on $\Pee^2_1$ and
$[u:v:w]$ on $\Pee^2_2$, we set
\begin{eqnarray} \label{p}
F\AND=\AND x^2yw^3+z^2yvw^2+z^3uw^2+(x+y)^3vw^2+p\ \ \in\ H^0(P,\OO(3,3))
\\ X\AND=\AND F^{-1}(0)\,\supset\,f\,=\,\{x^2y=0=u=v\}, \nonumber
\end{eqnarray}
for any $p$ vanishing on $f$ and sufficiently small (as above) for
$F^{-1}(0)$ to be smooth around $f$. The set of such $F$\,s form the
linear system
$$
H^0(P,\I_f(3,3))\subset H^0(P,\OO(3,3)),
$$
where $\I_f(3,3)$ denotes the sheaf of functions on $P$ vanishing on
$f$, twisted by the line bundle $\OO_P(3,3)$. It is easy to see that
the base locus of this linear system is just $f\subset P$: we want to
show that for any $x\in P\backslash f$, there exists an element of
the linear system not vanishing on $x$. This is clear if $x$
does not lie in the $\Pee^2_1$ fibre over $u=0=v$ that contains $f$,
and if it does we simply take one of the $F$\,s as above.

Thus, by Bertini's theorem [GH], the generic element of the linear system
has smooth zero locus away from $u=0=v$. In particular we may chose
$p$ arbitrarily small so that $X$ is both smooth on $u=0=v$, as above,
and away from it.

We now fix such a $p$, giving an $F$ and a smooth $X$, and study the
deformations of the line $C=\{x=0=u=v\}\subset f\subset X$ inside $X$.

\begin{Lemma}
$\nu_{C/X}\cong\OO_C(1)\oplus\OO_C(-3)$, so that $H^0(\nu_{C/X})\cong\C^2$. 
\end{Lemma}

\begin{Proof}
It is clear that $\nu_{C/P}=\nu_{C/\Pee^2_1}\oplus\nu_{\Pee^2_1/P}\cong
\OO_C(1)\oplus\OO_C\oplus\OO_C$, where $\Pee^2_1=\Pee^2_1\times\{[0:0:1]\}
\subset\Pee^2_1\times\Pee^2_2=P$. We set $w=1$ and work with
coordinates $[x:y:z],u,v$.

The standard exact sequence
$$
0\to\nu_{C/X}\to\nu_{C/P}\stackrel{dF\,}{\To}\OO(3,3)\res C\to0
$$
shows that $\nu_{C/X}$ is the kernel of the map
$$
\begin{CD}
\OO_C(1)\oplus\OO_C\oplus\OO_C @>\left(\left.{\partial F\over\partial
x}\right\res C,\left.{\partial F\over\partial u}\right\res C,\left.
{\partial F\over\partial v}\right\res C\right)\,>> \OO_C(3).
\end{CD}
$$
This map is $(0,z^3+\partial p/\partial u,z^2y+y^3+\partial
p/\partial v)$, since $x^2$ divides $p\res{u=0=v}$ ($p$ vanishes
scheme-theoretically on $f$) so that $\partial p/\partial x\res C=0$.
Thus the kernel is the first factor $\OO_C(1)$ plus
the kernel of the induced surjection $\OO_C\oplus\OO_C\to\OO_C(3)$,
which is $\OO_C(-3)$ by counting degrees.
\end{Proof}

Thus the first order deformations of $C$ are those inside the
$\Pee^2_1$ over $u=0=v$, through the normal bundle $\nu_{C/\Pee^2_1}
=\OO_C(1)$. Writing
$$
Y_\epsilon=Y\times\mathrm{Spec\,}\C[\epsilon]/(\epsilon^2),
$$
for the thickening of any space $Y$, the deformations give thickened
curves $C_\epsilon\into X_\epsilon$: that corresponding to the
section $s=\alpha y+\beta z\in H^0(\nu_{C/\Pee^2_1})$ is given by the
equations
$$
u=0=v,\quad x=\epsilon(\alpha y+\beta z).
$$
This lies inside $X_\epsilon$ since $x^2=\epsilon^2(\alpha
y+\beta z)^2$ and $\epsilon^2=0$. 

\begin{Lemma}
For $p$ in (\ref{p}) sufficiently small, and $C_\epsilon\into
X_\epsilon$ a curve defined by $s\in H^0(\nu_{C/\Pee^2_1})=
H^0(\nu_{C/X})$, $s$ is not in the image of the restriction map
$H^0(\nu_{C_\epsilon/X_\epsilon})\to H^0(\nu_{C/X})$ induced by
$\nu_{C_\epsilon/X_\epsilon}\res C=\nu_{C/X}$.
\end{Lemma}

\begin{Remarks}
Thus the curve $C_\epsilon$ cannot be deformed further: standard
deformation theory (see e.g. [Ka]) shows that the deformation
$C_\epsilon\into X_\epsilon$ given by
$s\in H^0(\nu_{C/X})$ can be extended to $C\times$\,Spec\,$\C[\epsilon]
/(\epsilon^3)\into X\times$\,Spec\,$\C[\epsilon]/(\epsilon^3)$ if and
only if $s$ is in the image of the above restriction map.

The following proof can be strengthened to show that in fact, for
$s\ne\alpha(y\pm iz)$ at least, $H^0(\nu_{C_\epsilon/X_\epsilon})=0$.
\end{Remarks}

\begin{Proof}
Repeating the above analysis of $\nu_{C/X}$ relative to the base
Spec\,$\C[\epsilon]/(\epsilon^2)$ (i.e. carrying $\epsilon$
as an extra variable, albeit one whose square is zero) shows that
$\nu_{C_\epsilon/X_\epsilon}$ is the kernel of the map
$$
\begin{CD}
\OO_{C_\epsilon}(1)\oplus\OO_{C_\epsilon}\oplus\OO_{C_\epsilon}
@>\left(\left.{\partial F\over\partial(x-\epsilon s)}
\right\res{C_\epsilon},\left.{\partial F\over\partial u}\right\res
{C_\epsilon},\left.{\partial F\over\partial v}\right\res
{C_\epsilon}\right)\,>> \OO_{C_\epsilon}(3),
\end{CD}
$$
where $s=\alpha y+\beta z$. Notice then that $\partial/\partial u$ and
$\partial/\partial v$ are taken holding $x-\epsilon(\alpha y+\beta
z)$, and not $x$, fixed. However $\epsilon^2=0$ simplifies things
to make the map
$$
\left(\epsilon(\alpha y+\beta z)(2y+q),z^3+\partial p/\partial u,
z^2y+3\epsilon(\alpha y+\beta z)y^2+y^3+\partial p/\partial v\right),
$$
with $q$ the polynomial $(\partial p/\partial x)/x$. We must show that
there are no sections of the kernel restricting on $\epsilon=0$ to
$s=(\alpha y+\beta z,0,0)$, i.e. that $(\alpha y+\beta z+\epsilon l,
\epsilon a,\epsilon b)$ is not in the kernel for all linear $l$ and
$a,\,b\in\C$.

Thus we must show there are no solutions in $a,\,b$ to
$$
2y(\alpha y+\beta z)^2+az^3+by(y^2+z^2)+\left[q(\alpha y+\beta z)^2
+a{\partial p\over\partial u}+b{\partial p\over\partial v}\right]=0,
$$
where the term in the square brackets can be made small with $p$.
Choose $p$ such that the coefficients of $z^3,\,z^2y,\,y^3$ in
$q,\,\partial p/\partial u$ and $\partial p/\partial v$ are less than
$1/10$. Then comparing coefficients of $z^3$ shows that $a=0$, while
the coefficients of $y(y^2+z^2)$ show that $b=0$ since $z^2+y^2$
cannot be written as a square $(\alpha y+\beta z)^2$.
\end{Proof}

\section{From curves to bundles}

Just as (codimension 1) divisors correspond to (rank 1) line
bundles, codimension 2 subschemes $Z$ of schemes $X$ sometimes
correspond to rank 2 bundles $E$ via zero sets of sections $s\in H^0(E)$.
If this is the case the subscheme will be a local complete intersection
(locally it is cut out by 2 independent equations given by the two
components of $s$ in a local trivialisation of $E$)
with normal bundle the restriction of $E$ to $Z$ (with the
isomorphism set up by the derivative of $s$ at $Z$). Letting
$L=\Lambda^2E$, wedging with $s$ gives maps $\OO\to E\to L$ giving an
exact sequence
\begin{equation} \label{ext}
0\to\OO\stackrel{s\,}{\to}E\to\I_Z\otimes L\to0,
\end{equation}
where $\I_Z$ is the ideal sheaf of functions vanishing on $Z$. We will
consider $Z$ connected for simplicity.

Thus, given $Z\subset X$, for it to be the zero locus of a section $s$
of a rank 2 bundle $E$, $Z$ must firstly be a local complete
intersection such that $\Lambda^2\nu_{Z/X}$ extends to a global line
bundle $L$ on $X$. If this is the case then to construct $E$ we want a
global extension (\ref{ext}) which is locally free (a vector
bundle). Away from $Z$ this is given by an element of $H^1(L^*)$,
which we may think of as an $L^*$-valued one-form that is
$\bar\partial$-closed away from $Z$. Local analysis on $Z$ shows that
the extension $E$ being locally free is equivalent to $\bar\partial$
of the form being a Dirac delta along $Z$. Thus the condition for the
global existence of such a vector bundle $E$ is that the class in
$H^2(L^*)$ defined by this Dirac delta is $\bar\partial$ of something,
i.e. that it is zero.

The upshot of all this is the Serre construction [GH], giving $(E,s)$
and an exact sequence (\ref{ext}) if and only if $[Z]\in H^2(L^*)$ is
zero. Here $[Z]\in H^2(L^*)$ is defined
as a linear functional on $H^{n-2}(L\otimes K_X)$ (the Serre dual of
$H^2(L^*)$, where $n=$\,dim\,$X$ and $K_X$ is the canonical bundle) by
restricting a form to $Z$ and integrating:
$$
H^{n-2}(L\otimes K_X)\to H^{n-2}_Z(\Lambda^2\nu_{Z/X}\otimes K_X)=
H^{n-2}_Z(K_Z)\stackrel{\int\,}{\to}\C.
$$

In particular if $H^2(L^*)=0$ then $E$ and $s$ exist.

Using this we now construct bundles $E\to X$ and $E_\epsilon\to
X_\epsilon$ from $C$ and $C_\epsilon$. We do the $C\subset X$ case,
but $C_\epsilon\subset X_\epsilon$ is identical.

$\Lambda^2\nu_{C/X}=\OO_C(-2)$ extends to all of $X$ as $L=\OO(-2,k)$
for any $k$; we choose $k=2$ to stand a chance of $E$ being
stable. The divisor exact sequence
\begin{equation} \label{div}
0\to\OO_P(k-3,l-3)\to\OO_P(k,l)\to\OO_X(k,l)\to0,
\end{equation}
and the cohomology of line bundles on $\Pee^2$ and hence $P$, show that
$H^2(L^*)=0$. Thus the Serre construction gives a rank 2 bundle $E$
and a section $s$ vanishing on $C$, and so the exact sequence
\begin{equation} \label{kos}
0\to\OO\stackrel{s\,}{\to}E\to\I_C(-2,2)\to0.
\end{equation}
Since the class $[C]\in H^4(X)$ is $\omega^2_2/3$, the Chern classes
of $E$ are $c_1=2\omega_2-2\omega_1$ and $c_2=\omega^2_2/3$. Thus the
Chern classes of $A=E(1,-1)$ are
\begin{eqnarray}
c_1(A) \AND=\AND 0, \nonumber \\
c_2(A) \AND=\AND \omega_1^2+\frac43\omega_2^2. \nonumber
\end{eqnarray}

\section{An obstructed bundle}

Repeated applications of the cohomology of the sequences $0\to\I_C\to\OO_X
\to\OO_C\to0$ and (\ref{div}), and Serre duality, give
$$
H^1(\OO_X(-2,2))=0 \so H^1(\I_C(-2,2))=0 \so H^1(E)=0=H^2(E^*),
$$ and
$$
H^1(\OO_X(2,-2))=0 \mathrm{\ and\ } H^1(\I_C)=0 \so H^1(E^*)=H^1(E(2,-2))=0.
$$
Thus tensoring (\ref{kos}) with $E^*=E(2,-2)$ shows that
$$
H^1(\End E)\cong H^1(E\otimes\I_C).
$$
(\ref{kos}) also shows that $H^0(E)$ is generated by $s$, so that the
restriction map $H^0(E)\to H^0(E\res C)$ is zero. Thus the sequence
$0\to E\otimes\I_C\to E\to E\res C\to0$ gives $H^1(E\otimes\I_C)\cong
H^0(E\res C)\cong H^0(\nu_{C/X})$ and we find that
\begin{eqnarray}
H^1(\End E) \AND\cong\AND H^0(\nu_{C/X}), \label{def1} \\ \label{def2}
H^1(\End E_\epsilon) \AND\cong\AND H^0(\nu_{C_\epsilon/X_\epsilon}).
\end{eqnarray}
(\ref{def1}) says that $E$ has a $\C^2\cong H^0(\nu_{C/X})$ of first
order deformations, while (\ref{def2}) says that for any such deformation
$E_\epsilon\to X_\epsilon$ given by some $s\in H^0(\nu_{C/X})$, $s$ is
not in the image of the restriction map $H^1(\End E_\epsilon)\to
H^1(\End E)$. 

Thus to prove Theorem \ref{thm} it is now sufficient to show the
slope-stability of $E$ with respect to the K\"ahler forms
$N\omega_1+\omega_2,\ N>1$. By replacing
subsheaves of $E$ by their double duals, it is enough to show that the
slopes of any line bundle subsheaves of $E$ are smaller than the slope
of $E$. The cohomology of $X$ shows that the only line bundles on $X$
are of the form $\OO(k,l)$, so it is sufficient to show that if
$E(k,l)$ has a section then deg\,$E(k,l)>0$.

The sequence (\ref{div}) shows that $\OO_X(a,b)$ has sections if and
only if $a\ge0$ and $b\ge0$, and for a section to vanish on $C$ we
must have $b\ge1$, so that $H^0(\I_C(a,b))\ne0$ if and only if
$a\ge0,\,b\ge1$. Thus we see from (\ref{kos}) that $E(k,l)$ cannot
have sections unless either $k,l\ge0$ or $k\ge2,\,l\ge-1$.

So it is sufficient to show that the degrees of $E$ and $E(2,-1)$ are
positive. They are
$$
c_1(E)\,.\,(N\omega_1+\omega_2)^2=6(N^2-1) \vspace{-3mm}
$$ and
$$
c_1(E(2,-1))\,.\,(N\omega_1+\omega_2)^2=12(2N+1),
$$
which are positive for $N>1$, as required.

\subsection*{Addendum}
Here we prove the moduli space is isomorphic to Spec\,$\C[\epsilon,\eta]/(\epsilon^2,\eta^2)$.
As Yukinobu Toda pointed out, 17 years later it is easier to prove the isomorphism between the deformation theory of $C$ and $E$ (and so $A=E(1,-1)$) by noting that $E$ is the inverse spherical twist of $\I_C(-2,2)$ about the spherical sheaf $\OO$. So we concentrate on describing ($C$'s component of) the moduli space of deformations of $C$ inside $X$.

Deformations of $C$ inside $P$ are lines described by the equations $x=ay+bz,\ u=U,\ v=V$ depending on the four parameters $a,b,U,V$. Restricting $X$'s defining equation $F$ (\ref{p}) to these lines gives
$$
(ay+bz)^2y+z^2Vy+z^3U+(x+y)^3V+p.
$$
Momentarily setting $p=0$ (which makes $X$ singular, but only away from $C$), this defines the section
$$
\Big(a^2+(1+a)^3V,\ 2ab+3(1+a)^2bV,\ b^2+V+3(1+a)b^2V,\ U+b^3V\Big)
$$
of the rank-4 bundle over $\C\!^4_{a,b,U,V}$ with basis the sections $y^3,y^2z,yz^2,z^3$ of $\OO(3)$ over any of the lines. The zero locus of these sections, localised near $a=b=U=V=0$, is our moduli space.

The first section gives $V=-a^2/(1+a)^3$ which we plug into the fourth $U=-b^3V$ to describe both in terms of $a,b$. Substituting into the third gives
\begin{equation}\label{fg}
b^2\ =\ -V(1+3(1+a)b^2)\ =\ a^2\frac{(1+3b^2+3ab^2)}{(1+a)^3}\,,
\end{equation}
and the second gives
$$
ab\left(2-\frac{3a}{1+a}\right)=0.
$$
Thus $ab=0$, which combined with (\ref{fg}) shows that all cubic terms in $a,b$ vanish. Therefore (\ref{fg}) simplifies to $b^2=a^2$. Changing variables to $\epsilon=a+ib,\,\eta=a-ib$ gives $\epsilon^2=0=\eta^2$ as claimed in the Introduction.

Thus when $p$ vanishes to order $\ge3$ on $C$ (which, by Bertini, is enough to ensure the smoothness of $X$) the moduli space is the same. For other small deformations $p$, the moduli space is contained in Spec\,$\C[\epsilon,\eta]/(\epsilon^2,\eta^2)$ by Theorem \ref{thm}, and can be no smaller since its length is deformation invariant and equals 4 when $p=0$.

\end{document}